\documentclass[12pt]{article}

\usepackage{amssymb,amsmath,exscale}

\usepackage{graphicx}

\usepackage{color}
\definecolor{marin}{rgb}   {0.,   0.3,   0.7} 
\definecolor{rouge}{rgb}   {0.8,   0.,   0.} 
\definecolor{sepia}{rgb}   {0.8,   0.5,   0.} 
\usepackage[colorlinks,citecolor=marin,linkcolor=rouge,
            bookmarksopen,
            bookmarksnumbered
           ]{hyperref}

\newtheorem{lemma}{Lemma}[section]

\newtheorem{theorem}[lemma]{Theorem}

\newtheorem{remark}[lemma]{Remark}
\newtheorem{example}[lemma]{Example}
\newtheorem{hypothesis}[lemma]{Hypothesis}
\newtheorem{notation}[lemma]{Notation}
\newtheorem{definition}[lemma]{Definition}
\newtheorem{conclusion}[lemma]{Conclusion}
\numberwithin{equation}{section}
\newcommand{\QED}{\mbox{}\hfill \raisebox{-0.2pt}{\rule{5.6pt}{6pt}\rule{0pt}{0pt}} 
          \medskip\par}

\newcommand{\dd}{\mathrm{d}}

\newcommand{\N}{\mathbb{N}}

\newcommand{\R}{\mathbb{R}}

\newcommand{\T}{\mathbb{T}}
\newcommand{\Z}{\mathbb{Z}}

\title{Quasi-periodic solutions of the 2D Euler equation}        
\author{Nicolas Crouseilles and Erwan Faou}       
\begin{document}
\maketitle
%\abstract{}

\abstract{We consider the two-dimensional Euler equation with periodic boundary conditions. We construct time quasi-periodic solutions of this equation made of localized travelling profiles with compact support propagating over a stationary state depending on only one variable. The direction of propagation is orthogonal to this variable, and the support is concentrated on flat strips of the stationary state. The frequencies of the solution are given by the locally constant velocities associated with the stationary state.    }

\section{Introduction}

We consider the two-dimensional Euler equation written in terms of vorticity
$$
\partial_t \omega + u \cdot \nabla \omega = 0,
$$
where $\omega(t,x,y) \in \R$, $\nabla = (\partial_x,\partial_y)^T$ with $(x,y) \in \T^2$ the two-dimensional torus $(\R \slash 2\pi\Z)^2$. The divergence free velocity field $u$ is given by the formula
\begin{equation}
\label{eq:J}
u = J \nabla \psi \quad \mbox{with} \quad \psi = \Delta^{-1} \omega, \quad \mbox{where}\quad 
J = \begin{pmatrix}
0 & - 1 \\ 1 & 0
\end{pmatrix}, 
\end{equation}
$J$ is the canonical symplectic matrix. Here $\Delta^{-1}$ is the inverse of the Laplace operator on functions with average $0$ on $\T^2$. 
We can rewrite this equation as 
\begin{equation}
\label{eq:euler}
\left\{
\begin{array}{l}
\partial_t \omega + \{ \psi, \omega\} = 0 ,\\[2ex]
\Delta \psi = \omega,
\end{array}
\right.
\end{equation}
with the 2D Poisson bracket for functions on $\T^2$: 
$$
\{f,g\} = (\partial_xf)(\partial_yg) - (\partial_yf)(\partial_xg). 
$$
%For two functionals $H(\omega)$ and $G(\omega)$, we set
%$$
%\{H,G\}_\omega = \int_{\T^2} \frac{\delta H}{\delta \omega}\left\{ \frac{\delta G}{\delta \omega}, \omega \right\} = -\{ G,H\}_\omega. 
%$$
%The Euler equation  \eqref{eq:euler} is a Hamiltonian PDE associated with this non canonical Poisson structure, and with Hamiltonian
%$$
%E(\omega) =   \frac{1}{(2\pi)^2}\int_{\T^2} \frac12 \Norm{u}{}^2 = -\frac{1}{2(2\pi)^2} \int_{\T^2} \omega (\Delta^{-1}\omega) = \Norm{\omega}{H^{-1}}^{2}, 
%$$
%which is quadratic in $\omega$. In other words, we can write \eqref{eq:euler} as
%$$
%\partial_t \omega + \left\{  \frac{\delta E}{\delta \omega}, \omega  \right\} = 0, 
%$$ 
%and from the definition of the Poisson structure, we observe that $E(\omega(t)) = E(\omega(0))$ for all time (preservation of the energy). Moreover, the flow is volume preserving in the sense that for all smooth functions $h:\R \mapsto \R$, we have 
%\begin{equation}
%\label{eq:casimirs}
%\forall\, t \geq 0\, \quad \int_{\T^2} h(\omega(t,x,y)) \dd x \dd y =  \int_{\T^2} h(\omega(0,x,y)) \dd x \dd y, 
%\end{equation}
%which expresses the preservation of the Casimirs of the Poisson structure. 
%
The equation \eqref{eq:euler} possesses many stationary states. For all functions $F:\R \mapsto \R$ and $\psi^0: \T^2 \to \R$ satisfying 
$
\Delta \psi^0 = F(\psi^0), 
$
then the couple of functions $\omega(t,x,y) = F(\psi^0(x,y))$ and $\psi(t,x,y) = \psi^0(x,y)$ solve \eqref{eq:euler}. We refer to \cite{Ambrosetti,CLMP92,CPR09,Hill} for further analysis of these particular solutions. 

Another class of stationary states are given by functions depending only on one variable (shear flows): for any smooth $V(x)$ periodic in $x$, the couple $\omega(t,x,y) = V''(x)$ and $\psi(t,x,y) = V(x)$ is solution of the 2D Euler equation. Note that $u(x) = (0,V'(x))^T$.

The goal of this paper is to construct solutions of \eqref{eq:euler} of the form 
\begin{equation}
\label{eq:armada}
\omega(t,x,y) = V''(x) + \sum_{k = 1}^K \Omega_k(x, y  -  v_k t), 
\end{equation}
where the functions $\Omega_k$ are localized around points $(x_k,y_k) \in \T^2$ such that $V''(x) = 0$ in a neighborhood of $x_k$ which correspond to flat strips of the stationary state $V''(x)$. The points $y_k$ are arbitrary points in $[0,2\pi]$, and $v_k = V'(x_k)$ is the locally constant velocity associated with $V''(x)$. The profiles $\Omega_k$ are constructed as stationary states of the 2D Euler on $\R^2$ with radial symmetry around $(x_k,y_k)$, compact support and zero average. 

This very simple and explicit construction allows to construct quasi-periodic solutions to the 2D Euler equation corresponding to invariant tori of any given dimension in the dynamics. 
Let us stress that this is in surprising contrast with the traditional situation in nonlinear Hamiltonian PDEs such as Schr\"odinger and wave equations for which the construction is much more difficult and requires in general the use of Nash-Moser or KAM iterations, see for instance \cite{Craig93,Bourgain,poschel,EK,Berti}. 
%\begin{remark}
%These submarines quasi-solutions have been detected by numerical simulations. An extensive numerical study of the long time behavior of numerical schemes applied to the 2D periodic Euler equations, as well as the analysis of the quasi-periodic nature of submarines-like solutions is currently under study by the authors, in collaboration with Sergei Kuksin. 
%\end{remark}
%

%#############################
%#############################

\section{Construction}
Let $V(x)$ be a periodic function with zero average in $x \in [0,2\pi]$. We make the following assumptions on $V$:

\begin{hypothesis}
\label{hyp}
Let $K \in \N$. For all $k = 1,\ldots,K$, there exist $a_k < b_k$ in $[0,2\pi]$ such that for all $j,k \in \{ 1,\ldots,K\}$, we have $[a_j,b_j] \cap [a_k,b_k] = \emptyset$, and such that $V''(x) = 0$ for $x \in [a_k,b_k]$. 
\end{hypothesis}

This hypothesis implies that for $x \in [a_k,b_k]$, $V'(x) =: v_k$ is constant. 
Let us seek a solution $\omega(t,x,y)$ under the form 
$$
\omega(t,x,y) = V''(x) + \sum_{k = 1}^K \chi_k(t,x,y), 
$$
where for all $k$, $\chi_k(t,x,y)$ is of zero average, and the support of $\chi_k$ and $\Delta^{-1} \chi_k$ is included in $]a_k,b_k[$. This implies in particular that for all $j$ and $k$, 
$$
\{ \Delta^{-1}\chi_j,\chi_k\} = 0 \quad \mbox{for}\quad j \neq k. 
$$
Inserting this decomposition into \eqref{eq:euler}, we thus obtain 
\begin{equation}
\label{eq:lin}
\sum_{k = 1}^K \Big(\partial_t \chi_k + \{ \Delta^{-1} V'',\chi_k\} + \{ \Delta^{-1} \chi_k, V''\} + \{ \Delta^{-1} \chi_k ,\chi_k\}\Big) = 0. 
\end{equation}
%%########
%A PARTIR DE LA
%%########

%Let us write $\omega(t,x,y) = V''(x) + \chi(t,x,y)$. The equation for $\chi(t,x,y)$ is given by 
%$$
%\dot \chi + \{ \Delta^{-1} V'',\chi\} + \{ \Delta^{-1} \chi, V''\} + \{ \Delta^{-1} \chi ,\chi\} = 0, 
%$$
%or equivalently
%$$
%\dot \chi + V'(x) \partial_y \chi - V'''(x) \partial_y \Delta^{-1}\chi + \{ \Delta^{-1} \chi ,\chi\} = 0. 
%$$
We seek for travelling wave solutions $\chi_k(t,x,y)  = \Omega_k(x,y - v_k t)$. Equation (\ref{eq:lin}) then becomes 
$$
\sum_{k = 1}^K \Big((V'(x) - v_k)  \partial_y \Omega_k - V'''(x) \partial_y \Delta^{-1}\Omega_k + \{ \Delta^{-1} \Omega_k ,\Omega_k\} \Big)= 0, 
$$
where $\Omega_k(x,y)$ and $\Psi_k(x,y) = \Delta^{-1} \Omega_k$ have compact support in $]a_k,b_k[ \times \T$. 
Using Hypothesis \ref{hyp}, we have: $V'''(x) = 0$ and $V'(x)=v_k$ for $x \in [a_k,b_k] \subset [0,2\pi]$. 
As the intervals $[a_k,b_k]$ are pairwise disjoints, the system of equations to solve is hence: For all $k = 1,\ldots,K$,  
\begin{equation}
\label{eq:stat}
 \{ \Psi_k,\Omega_k\}(x, y) = 0 \quad\mbox{and}\quad \Delta \Psi_k(x,y) = \Omega_k(x,y), \quad (x, y) \in [a_k,b_k] \times \T. 
\end{equation}
In other words, the couple $(\Omega_k,\Psi_k)$ is a smooth stationary state of the Euler equation with support on the flat strip $]a_k,b_k[ \times \T \subset \T^2$.

 To prove the existence of such functions $(\Omega_k,\Psi_k)$, take $k \in \{1,\ldots,K\}$, and fix 
$(x_k,y_k) \in  \, ]a_k,b_k[ \times [0,2\pi]$.  
Let us perform the local action-angle change of coordinates  $x - x_k = \sqrt{2r}\cos\theta$, $y -y_k = \sqrt{2r} \sin\theta$ for $(x,y)$ close enough to $(x_k,y_k)$. The jacobian matrix is
$$
M = \begin{pmatrix}
 \frac{1}{\sqrt{2r}} \cos \theta & \frac{1}{\sqrt{2r}} \sin \theta \\
-\sqrt{2r}\sin \theta & \sqrt{2r} \cos \theta
\end{pmatrix}, 
$$
and we verify that this matrix satisfies $M^T J M = J$, 
%$J=
%\begin{pmatrix}
%0 & -1\\ 
%1 & 0 
%\end{pmatrix},
%$ 
%the change of variables satisfies 
%$$
%M^T J M = 
% \begin{pmatrix}
%\frac{1}{\sqrt{2r}} \cos \theta & -\sqrt{2r}\sin \theta \\
%\frac{1}{\sqrt{2r}} \sin \theta  & \sqrt{2r} \cos \theta
%\end{pmatrix}
%\begin{pmatrix}
%\sqrt{2r}\sin \theta & -\sqrt{2r} \cos \theta\\
%\frac{1}{\sqrt{2r}} \cos \theta & \frac{1}{\sqrt{2r}} \sin \theta 
%\end{pmatrix}
% = J. 
%$$
which means that the change of coordinate is symplectic. Hence this transformation preserves the Poisson bracket and if we define $\tilde\Psi_k(r,\theta) = \Psi_k(x,y)$ and  $\tilde\Omega_k(r,\theta) = \Omega_k(x,y)$ we have 
\begin{align*}
\{ \Psi_k,\Omega_k\} &=(\partial_x \Psi_k )(\partial_y \Omega_k )- (\partial_y \Psi_k)( \partial_x \Omega_k) \\
&=(\partial_r\tilde \Psi_k )(\partial_\theta \tilde\Omega_k) - (\partial_\theta \tilde\Psi_k)( \partial_r \tilde\Omega_k) = \{ \tilde \Psi_k,\tilde \Omega_k \}. 
\end{align*}
%In coordinate $(r,\theta)$, the metric is given by 
%$$
%g_{rr} = \frac1{2r}, \quad g_{r\theta}= 0, \quad g_{\theta \theta} = 2r, 
%$$
%and the determinant $g$ is equal to $1$. 
The Laplace operator in coordinate $(r,\theta)$ is 
$$
\Delta f = 2(\partial_r f +r\partial_{rr} f) +  \frac{1}{2r} \partial_{\theta\theta} f.
$$
%$$
%\Delta \psi = \frac{1}{\sqrt{g}} \partial_i (\sqrt{g} g^{ij} \partial_j \psi), 
%$$
%is [{\bf QUESTION : Si $g=M M^T$ tel que $(M M^T)_{i,j} = g_{i,j}$, 
%on a $g^{i j} = ((M^T M)^{-1})_{i,j}$}. Donc $g^{rr}= 2r, g^{\theta\theta}=1/(2r)$. ] 
%%\begin{eqnarray*}
%%\Delta \psi  &=& \partial_r (g^{rr} \partial_r \psi) +  \partial_\theta ( g^{\theta\theta} \partial_\theta \psi)\\
%%&= & \partial_r (\sqrt{2r}\partial_r \psi) +  \partial_\theta ( \frac{1}{2\sqrt{r}} \partial_\theta \psi)\\
%%{\bf rajout} &= & \frac{1}{\sqrt{2r}}\partial_r \psi + \sqrt{2r}\partial_{rr}\psi+ \frac{1}{2\sqrt{r}}  \partial_{\theta\theta} \psi\\
%%{\bf manque \, 1/\sqrt{2} \, devant \, \partial_{\theta\theta}}
%%&=& \frac{1}{\sqrt{2r}} (\partial_{r} + \partial_{\theta \theta})\psi + \sqrt{2r} \partial_{rr}\psi
%%\end{eqnarray*}
%
%{\bf CORRECTION} 
%\begin{eqnarray*}
%\Delta \psi  &=& \partial_r (g^{rr} \partial_r \psi) +  \partial_\theta ( g^{\theta\theta} \partial_\theta \psi)\\
% &=& \partial_r (2r \partial_r \psi) +  \partial_\theta \Big( \frac{1}{2r}\partial_\theta \psi\Big)\\
%  &=& 2(\partial_r \psi +r\partial_{rr} \psi) +  \frac{1}{2r} \partial_{\theta\theta} \psi\\
%\end{eqnarray*}
%
%\medskip
%\medskip

Now take a smooth function $\tilde \Psi_k(r)$ on $[0,+\infty]$ such that $\int_{0}^{+\infty} \tilde \Psi_k(r) \dd r= 0$ 
and with compact support in 
$r < \min(|x_k-a_k|, |x_k-b_k|)$.  
Set 
$$
%\eta_k(r,\theta) = \frac{1}{\sqrt{2r}} \psi_k'(r)  + \sqrt{2r} \psi_k''(r) = \Delta \psi_k, \quad \quad k=1, \dots, K. 
\tilde \Omega_k(r) = \Delta\tilde  \Psi_k(r) =  2( \tilde \Psi_k'(r)  + r \tilde \Psi_k''(r)). 
$$
We verify that $\tilde \Omega_k(r)$ has compact support, that $\int_0^\infty \tilde \Omega_k(r) \dd r = 0$, and that $\{ \tilde \Psi_k,\tilde \Omega_k\} = 0$ as $\tilde\Psi_k$ and $\tilde\Omega_k$ only depend on $r$. 
Going back to the variables $(x,y)$, we verify that the function $(\Omega_k(x,y),\Psi_k(x,y))$ extended by $0$ outside the strip $]a_k,b_k[ \times \T$ satisfies \eqref{eq:stat} and are of zero average on $\T^2$. 

This proves the following result: 
%We then have [{\bf $\eta_k$ only depends on $r$ and }] as $\eta_k$ is of zero average  [{\bf peut-on enlever la dependance en $\theta$ de $\eta_k$ ??}]
%\begin{equation}
%\label{eq:laprt}
%\Delta^{-1} \eta_k (r,\theta) = \psi_k(r),  \quad k=1, \dots, K, 
%\end{equation}
%and hence $\{\Delta^{-1} \eta_k, \eta_k\} = \{ \psi_k(r), \eta_k(r) \}=  0$. 
%
%In the end, we have constructed a solution 
%$$
%\omega(x,y) = V''(x) + \sum_{k=1}^K \eta_k( x-x_k, y - y_k - v_k t), 
%$$ 
%to the 2D Euler equation, with $\eta$ with support included in $|x - x_k| \leq \sigma$. 
%As $x_k \in ]a_k,b_k[$, this implies $x \in [a_k,b_k]$ is $\sigma$ is sufficiently small 
%(typically $0<\sigma < \min(|x_k-a_k|, |x_k-b_k|)$), and $x$ belongs to an interval 
%on which $V'(x)$ is constant. 

\begin{theorem}
Assume that $V=V(x)$ is a periodic function with zero average satisfying Hypothesis \ref{hyp}. Then for all $k = 1,\ldots,K$, there exists points $(x_k,y_k) \in ]a_k,b_k[ \times \T$ and zero average functions $(\Omega_k, \Psi_k)$ with compact support in $]a_k,b_k[ \times \T$ and radial symmetry around $(x_k,y_k)$  such that the couple
$$
\omega(t,x,y) = V''(x) + \sum_{k=1}^K \Omega_k( x, y  - v_k t), 
$$ 
and 
$$
\psi(t,x, y) = V(x) +  \sum_{k=1}^K \Psi_k( x, y  - v_k t), 
$$
with $v_k=V'(x_k)$ is a quasi-periodic solution of (\ref{eq:euler}). 
\end{theorem}

\paragraph{Acknowledgements:}
It is a great pleasure to thank Sergei Kuksin and Nicolas Depauw for many helpful discussions.

\noindent {\bf Authors address:} 

\vskip 2ex

\noindent N. Crouseilles and E. Faou, INRIA and ENS Cachan Bretagne, Avenue Robert Schumann, 35170 Bruz, France.

\noindent {\tt Nicolas.Crouseilles@inria.fr}, 
{\tt Erwan.Faou@inria.fr}

\end{document}